\documentclass[12pt]{amsart}
\usepackage{amsmath,amsthm,epsfig,amssymb}
\usepackage{amsthm,amsfonts,amssymb,amscd}
\usepackage[all]{xy}
\newtheorem{thm}[subsection]{Theorem}

\newtheorem{lem}[subsection]{Lemma}
\newtheorem{corol}[subsection]{Corollary}
\newtheorem{rem}[subsection]{Remark}
\theoremstyle{definition}
\newtheorem{Def}[subsection]{Definition}

\newtheorem{exam}[subsection]{Example}
\newtheorem{proposition-definition}[subsection]{Proposition-Definition}

\newcommand{\CC}{{\mathbb C}}

\newcommand{\ZZ}{{\mathbb Z}}

\newcommand{\PP}{{\mathbb P}}

\newcommand{\OOO}{{\mathcal O}}

\author{A. El Mazouni}
\address{ Univ. Artois, UR 2462, Laboratoire de 
Mathématiques de Lens (LML), F-62300 Lens, France.}
\email{abdelghani.elmazouni@univ-artois.fr}

\author{D. S. Nagaraj}
\address{Indian Institute of Science Education and Research, 
Rami Reddy Nagar, Karakambadi Road,
Mangalam (P.O.), Tirupati - 517507,
Andhra Pradesh, INDIA.}
\email{dsn@iisertirupati.ac.in}
\title{Hyperplane sections of the projective bundle associated to the tangent bundle of $\mathbb{P}^2.$  }
\date{}
\begin{document}

\begin{abstract}  
The aim of the note is to give a complete description of all the hyperplane sections of the projective bundle associated 
to the tangent bundle of ${\PP}^2$ under its natural embedding in $\mathbb{P}^7.$ As an application one obtains 
a description of all possible deformation in $\text{SL}_3/B$  of the co-dimension one subscheme which is the union
of two fundamental Schubert divisors. 

\end{abstract}

\maketitle

\section{Introduction} \setcounter{page}{1} 
Throughout we work over the field of complex numbers $\mathbb{C}.$
Let $X$ be a projective variety embedded in a projective space  $\mathbb{P}^n.$ 
Recently several researchers are lead to study the hyperplane sections of
$X$ in $\mathbb{P}^n$(see \cite{BMSS},\cite{CHMN}, \cite{GIV}). 
Given $X \subset \mathbb{P}^n,$ an interesting and challenging 
problem is to classify  the subschemes $H \cap X$  of $X$ 
for varying hyperplanes $H \subset \mathbb{P}^n.$ 
This problem is in general very difficult. 
For example, if we denote by $V(2, d)$ the image of $\mathbb{P}^2$ in its 
$d$-tuple embedding, then 
one has the description for all the hyperplane sections only for $d \leq 3.$
For $d = 1, 2$ it is easy and for $d = 3,$ see \cite{Ne}. In the case $d=2$ there are
three distinct classes of subschemes  and in the case $d=3$  there are
eight distinct classes of subschemes are there.

Let $\mathbb{P}(T_{\mathbb{P}^2})$ be the projective bundle
associated to the tangent bundle of $T_{\mathbb{P}^2}.$ This threefold is 
naturally embedded in 
$\mathbb{P}^7.$  The aim of this short note is to describe explicitly all 
possible types of hyperplane 
sections of this embedding. Existence of such a complete description is  very rare.

We have the following:

\begin{thm}\label{thm1}
A hyperplane section of $\mathbb{P}(T_{\mathbb{P}^2})$ in its natural 
embedding in $\PP^7$ is either 

 i)  an irreducible surface of degree six in $\PP^6.$  
 
 or
 
 ii) the union of two non-singular degree three surfaces.
 
 If the hyperplane section is irreducible then it is 
 either a non-singular Del Pezzo surface or a singular rational normal surface.
 Moreover, if it  is a singular rational normal surface then it is non-singular except one point
 at which  the singularity is either  $A_1$ type (i.e., quadratic cone singularity
 with the  local equation $Y^2+XZ= 0$ ) or $A_2$ type (i.e.,singularity
 with the local equation $Y^3+XZ =0$).
 
 If the hyperplane section is union of two non-singular degree three surfaces then
 their intersection is either one irreducible non-singular rational curve or a union of 
 two non-singular rational curves intersecting at a point. 
 
\end{thm}

Theorem \ref{thm1} gives us all possible deformations of the codimension one subscheme 
$H$ which is union of two fundamental co-dimension one Schubert varieties in
the homogenous space ${\text{SL}_3(\mathbb{C})}/B,$ where $\text{SL}_3(\mathbb{C})$ is the group of 
all $3\times 3$  matrices of determinant one over the field of 
complex numbers and $B$ is the subgroup of upper triangular matrices (see, Remark \ref{deform}).

\section{Proof of the Theorem}
Throughout this note we use standard notations, see for example \cite{GH} or \cite{Ha}. 

\begin{Def} For a subscheme $Z,$ of dimension zero, of the projective  space $\mathbb{P}^n,$  the $\CC$ vector space
dimension of the $\CC$ algebra $\OOO_Z$  is called the  length 
of the subscheme $Z$ and is denoted by $\ell(Z).$ i.e.,
$$\ell(Z) := \text{dim}_{\CC} \,\OOO_Z.$$  
\end{Def}
The proof of Theorem \ref{thm1}  depends on the following:

\begin{thm}\label{thm2}
Let $s$ be a non-zero section of the tangent bundle $T_{\mathbb{P}^2}$ of 
${\mathbb{P}^2}.$ Then the  scheme $V$ of zeros of $s$ is of one of the following types:

 a)  $V$ is reduced and consists of three distinct points,
 
 b) $V$ is non reduced and is  a union of a reduced point and a non reduced subscheme of length two,
 
 c) $V$ is non reduced subscheme of length three, supported on  a single point,
 
 d) $V$ is the union of a line and a point not on the line,
 
 e) $V$ is a line with an embedded point on it.
\end{thm} 

 {\bf Proof}: Let $\{X, Y, Z \}$ be a basis of the space of  sections  of the 
 line bundle $ {\mathcal O}_{\mathbb{P}^2}(1) .$  We consider $X, Y$ and $Z$ as the
 variables  giving the homogeneous 
 coordinates on $\mathbb{P}^2.$ On $\mathbb{P}^2$ we have the Euler sequence
 \begin{equation}\label{euler}
 0 \to {\mathcal O}_{\mathbb{P}^2} \to {\mathcal O}_{\mathbb{P}^2}(1)^3 
 \to T_{\mathbb{P}^2} \to 0, 
 \end{equation}
  where the map ${\mathcal O}_{\mathbb{P}^2} \to {\mathcal O}_{\mathbb{P}^2}(1)^3 $ 
  is given by
 $1 \mapsto (X,Y,Z).$   Any section of 
 ${\mathcal O}_{\mathbb{P}^2}(1)^3 $ gives a section of $T_{\mathbb{P}^2} $ 
 by projection. 
 Since  the group
 $$H^1({\mathcal O}_{\mathbb{P}^2} )= 0$$ 
 every section of  $T_{\mathbb{P}^2} $  is the  image of a section of 
 ${\mathcal O}_{\mathbb{P}^2}(1)^3. $
 A section of
 ${\mathcal O}_{\mathbb{P}^2}(1)^3 $ is an ordered triple $(f_1, f_2,f_3)$ 
 where $f_i \,(1\leq i\leq 3)$ are 
 linear forms in the variables $X,Y,Z.$ Two sections $(f_1, f_2,f_3)$ and 
 $(g_1, g_2,g_3)$ of 
 ${\mathcal O}_{\mathbb{P}^2}(1)^3 $  map to the same section of $T_{\mathbb{P}^2} $ if
 the difference of these two sections is a scalar multiple of $(X,Y,Z).$ Let $s$ be a 
 non-zero section of $T_{\mathbb{P}^2}$ and $(f_1, f_2,f_3)$ be a section of 
 ${\mathcal O}_{\mathbb{P}^2}(1)^3$
 which maps to $s$ under the surjection given by the Euler sequence 
 (\ref{euler}). Clearly  the scheme $V$ of zeros of $s$
 is equal to the subscheme of $\mathbb{P}^2$ on which the  two sections 
 $(X, Y, Z)$  and  $(f_1, f_2,f_3)$ of 
 ${\mathcal O}_{\mathbb{P}^2}(1)^3 $ are dependent, namely the scheme defined by the vanishing of
 the two by two minors of the matrix
 \[
 \left( \begin{array}{ccc}
 X&Y&Z \\
 f_1&f_2 &f_3
 \end{array} \right)
 \]
 i.e., the scheme $V$ defined by the common  zeros  of the polynomials 
 \begin{equation}\label{eq0}
 Xf_2-Yf_1, Xf_3-Zf_1,Yf_3-Zf_2.
 \end{equation}
 Since the second Chern class $c_2({\mathcal O}_{\mathbb{P}^2}(1)^3) = 3[R],$ 
 where $[R]$ is the cohomology class of a point $R \in \mathbb{P}^2,$
 we see that if the support of $V$ 
 is a finite set, then it must be one of the types  a), b) or c). Moreover as $s$ is 
 non zero we see that $V$ is defined 
 by at least  two linearly independent quadrics and hence it cannot 
 contain a conic as a subscheme.
  On the other hand if the restrictions of the sections
 $(f_1, f_2,f_3)$ and $(X,Y,Z)$ to a line $\ell$  are linearly dependent, then 
 the line $\ell$ is contained in $ V.$ In this 
 case we see  that $V$ is of the type $d)$ or $e)$ as at least two of the 
 quadrics in (\ref{eq0}) are
 linearly independent. $\hfill{\Box}$
 
 \begin{corol}\label{cor1}
 If $s$ is a non-zero holomorphic vector field on the projective plane $\mathbb{P}^2$ which vanishes on a
 positive dimensional subscheme  then that subscheme  is necessarily a  line i.e., given by zeros of a homogenous polynomial
 of degree one in three variables. 
 More over, if $s$ vanishes on a line then it can vanish at most one another point on $\mathbb{P}^2.$
 \end{corol}
 
 {\bf Proof}: Since a holomorphic section of the Tangent bundle is a holomorphic vector field 
 corollary is an immediate consequence of \ref{thm2}. $\hfill{\Box}$
  
  \exam {Here we give examples to show that  all the types 
  mentioned in Theorem \ref{thm2} do occur. 
  For $\lambda, \mu \in \mathbb{C},$ consider the section $( X,  \lambda Y, \mu Z)$ of 
  $ {\mathcal O}_{\mathbb{P}^2}(1)^3$ and denote by $s$
  the induced section of  $T_{\mathbb{P}^2} $ under the  surjection
  \[ {\rm H}^0({\mathcal O}_{\mathbb{P}^2}(1)^3 )\to {\rm H}^0(T_{\mathbb{P}^2} )\to 0 \]
  obtained by the exact sequence (\ref{euler}). Note that $\lambda = \mu =1$ 
  if and only if $s$ corresponds to the zero section of 
  $T_{\mathbb{P}^2}.$ The scheme $V$ of zeros of the section $s$ is given by
  the vanishing of the quadrics
  \begin{equation}\label{eq1}
 (\lambda-1)XY, (\mu -1)XZ, (\mu -\lambda)YZ.
 \end{equation}
 
 If $\lambda -1, \mu -1, ( \mu -\lambda) $ are all non zero, then 
 $$V=\{(1,0,0),(0,1,0), (0,0,1)\}$$ 
 consists of three distinct  points and hence is of type $a)$ of 
 Theorem \ref{thm2}.  If $\lambda = 0$ and $\mu = 1$ then $V$ is the union of the line $Y= 0$ 
 and the point $(0,1,0)$ and the 
 scheme is of type $d)$ of 
 Theorem \ref{thm2}. 
 
 If $s$ is the section defined by $(Y, 0,0)$ then the zero scheme of $s$ 
 is given by the 
 vanishing of the quadrics $Y^2$ and $YZ$ which is of  type $e)$ of 
 Theorem \ref{thm2}.
 
 If $s$ is the section defined by $(X, Z,0)$ then the zero scheme of $s$ 
 is given by the 
 vanishing of the quadrics $XZ-YX, XZ$ and $Z^2$ which is of  
 type $b)$ of Theorem \ref{thm2}.

If $s$ is the section defined by $(Y, Z,0)$ then the zero scheme of $s$ 
is given by the 
 vanishing of the quadrics $XZ-Y^2, YZ$ and $Z^2$ which is of  
 type $c)$ of Theorem \ref{thm2}.
 $\hfill{\Box}$

 \vskip3mm

Let $E$ be a vector bundle of rank two on a smooth projective surface $X.$ Let
 $\mathbb{P}(E)$ be the projective bundle
associated to the bundle  $E$  and $\pi : \mathbb{P}(E) \to X$ 
be the natural projection.
Let $ {\mathcal O}_{\mathbb{P}(E)}(1)$ be the 
relative ample line bundle quotient of $\pi^*(E)$ and let 
$\phi$ be the natural isomorphism 
\[ {\rm H}^0(X, E) \simeq 
{\rm H}^0({\mathbb{P}(E)}, {\mathcal O}_{\mathbb{P}(E)}(1)) .\]

 \begin{lem}\label{lem1}
 For a non zero section  $\sigma$ of  $E,$  let $V_0(\sigma)$  (resp.  $V_1(\sigma)$)
 be
 the union of zero dimensional  (resp. one dimensional) components of the 
 subscheme of $X$ defined by zeros of 
 $\sigma.$ Let $s\neq 0$ be a section of $E$ and $\phi(s)$ be the corresponding section 
 of ${\mathcal O}_{\mathbb{P}(E)}(1). $ 
 The scheme of zeros of a section $\phi(s)$ is equal to $W_0\cup W_1,$
where $W_0$ (resp.$ W_1$) is a subscheme of $\mathbb{P}(E)$  
 isomorphic to the blow-up of $X$ 
along $V_0(s)$ (resp. is isomorphic to ${\mathbb{P}(E|_{V_1(s)}}))$.
 \end{lem} 
  
  {\bf Proof:} To prove the  Lemma,  first we observe: 
  
  i)  If the section $s$ is non zero at a point $x$ 
  the subspace 
  generated by it  defines a unique point in the projective line 
  $\mathbb{P}(E_x)$ at which the section 
  $\phi(s)$ vanishes.
  
  ii)  Let $U$ be an affine open subset of $X$ such that $ E|_U \simeq \OOO_U^2.$  
  Now $s|_U$ gives a section of $\OOO_U^2,$ say $(f_1,f_2),$
  then the associated section $\phi(s)|_U$ of 
  ${\mathcal O}_{\PP(\OOO_U^2)}(1)$ is $f_1X_1+f_2X_2,$ where 
  $X_1,X_2\in {\text{H}}^0(\OOO_{\PP^1}(1))$ is a basis. 
  Thus it follows that the subscheme of vanishing of $\phi(s)|_U$ 
  is the same as that of $f_1X_1+f_2X_2.$ 
  
   If $I$ is the ideal generated by $f_1,f_2$ in ${\rm H}^0(\OOO_U),$ then
  $$I = \cap_{i=1}^rQ_i \cap_{i=1}^s F_i,$$
    is the primary decomposition  of $I,$ where $Q_i , 1\leq i\leq r$ (resp. $F_i, 1\leq i \leq s$)
    are primary ideals corresponding to irreducible components of $V_0(s)\cap U$ (resp. $V_1(s) \cap U$). 
    Since $U$ is non-singular, the components of $V_1(s)$ are locally principal divisors. That is, the subschemes
    $$Z_i=V(F_i) ,  1\leq i \leq s, $$  
    of $U$  are  locally principal. Hence both $f_1$ and $f_2$ vanish along the closed sub-scheme 
    $Z_i \subset V_1(s)\cap U, 1\leq i \leq s,$ and we see that 
    $${\mathbb{P}(E|_{V_1(s)}})\cap \pi^{-1}(U) = W_1 \cap \pi^{-1}(U),$$
  where $\pi: \mathbb{P}(E) \to X$ is the projection map.  By i) we see that the blow-up of $U$ along the subscheme
  $V_0(s)\cap U$  is equal to $W_0 \cap \pi^{-1}(U).$
  
  Since $X$ is covered by open sets of the form $U$ in ii) the Lemma follows immediately. 
  
 \vskip3mm
 
  \begin{lem}\label{lem2}
 Let $s$ be a non-zero section of the tangent bundle $T_{\mathbb{P}^2}$ of 
${\mathbb{P}^2}.$ Assume that the support of the scheme $ V (\subset {\mathbb{P}^2})$ of zeros of $s$ is at most 
two points.  Then either 

i) $V$ is a union of a reduced scheme supported at a point and a non reduced subscheme $V_1$ of length two
with $\OOO_{V_1}$  isomorphic to $\CC[X,Y]/{(X, Y^2)},$

or

ii) $V$ is a non reduced scheme of length three supported at a point and 
$\OOO_V$ isomorphic to $\CC[X,Y]/{(X, Y^3)}.$
 
 \end{lem}
 
 {\bf Proof}: Notations as in the proof of  Theorem \ref{thm2}. If $(f_1, f_2, f_3)$ is a section of 
 ${\mathcal O}_{\mathbb{P}^2}(1)^3 $  that map to the section $s$ of $T_{\mathbb{P}^2} $
 then  
 $$ V  = V(f_1Y-f_2X, f_1Z-f_3X, f_2Z-f_3Y).$$
 The assumption on the support of  the subscheme  $V,$   implies that 
 $V$ has to be either of type b) or c) of the Theorem \ref{thm2}. Observe that the cohomology class  of the zero 
 dimensional scheme $V$ is $3[R],$ where $[R] \in H^4(\mathbb{P}^2, \ZZ)$ is the class of a point 
 $ R \in \mathbb{P}^2.$ Thus length of $\CC$ algebra $\OOO_V$ is three (i.e., 
 $\text{dim}_{\CC}{\OOO_{V} }= 3$).
 If  $V$ is of the type b) of  Theorem \ref{thm2}, then
  $$\OOO_V = \OOO_Q/{m_Q} \oplus \OOO_P/I,$$
    where
  $Q, P \in \mathbb{P}^2 $ are two distinct points  and $m_Q$ is the maximal ideal of $Q$
  and $I$ is an ideal  contained in the maximal ideal $m_P$  of $Q$ satisfying 
  $$\text{dim}_{\CC}{\OOO_{P}/{I}} = 2.$$
  It is easy to see that any subschemes of length two supported on single point of $\mathbb{P}^2$ is isomorphic to 
  $\text{Spec}(\CC[X,Y]/{(X, Y^2)}),$ we conclude that if $V$ is of type b) of the Theorem \ref{thm2} then it
  corresponds to  i) of the Lemma.
 
 If $V$ is of the type c) of  Theorem \ref{thm2}, then
  $$\OOO_V =  \OOO_P/I,$$
    where
  $P \in \mathbb{P}^2 $ is a point 
  and $I$ is an ideal  contained in the maximal ideal $m_P$ of $P$ satisfying 
  $$\text{dim}_{\CC}{\OOO_{P}/{I}} = 3.$$
  Any subscheme whose support is a single point of $\mathbb{P}^2$ and is of length three is isomorphic to either
  $$\text{Spec}(\CC[X,Y]/{(X, Y^3)})$$
  or 
  $$\text{Spec}(\CC[X,Y]/{(X^2, Y^2, XY)}).$$
   Claim: If $V$ is of the form c) of Theorem \ref{thm2} then it
  corresponds to  ii) of the Lemma.
  To prove the claim, with out loss of generality, we assume that $P=(0,0,1).$ 
  Then the subscheme $V $ is contained in the affine open set $ U= \{ Z \neq 0\}$ of
  $\mathbb{P}^2.$ and  $ V  = V(g_1Y-g_2X, g_1-g_3X, g_2-g_3Y) $ where 
  $g_i, (1\leq i\leq 3)$ are homogenous polynomials in $X,Y$ of degree less or equal to one
  and the ideal $I = (g_1Y-g_2X, g_1-g_3X, g_2-g_3Y).$ By using  assumptions that 
  $P =(0,0,1)$ and $V$ is of the form  c) of Theorem \ref{thm2}  we conclude
  $\OOO_V  \simeq \CC[X,Y]/{(X, Y^3)}.$ This proves the claim. 
  
  This completes the proof of the Lemma. $\hfill{\Box}$
  
\vskip3mm

\begin{lem}\label{lem3}
 Let $s$ be a non-zero section of the tangent bundle $T_{\mathbb{P}^2}$ of 
${\mathbb{P}^2},$ and $\phi(s)$ be the corresponding section of the line bundle
$ {\mathcal O}_{\mathbb{P}(T_{\mathbb{P}^2})}(1).$
Assume that the support of the scheme $ V (\subset {\mathbb{P}^2})$ of zeros of $s$ is at most 
two points.   

i) If $V$ is a union of a reduced point $Q$ and a non reduced scheme of length 
two supported  at point $P \in \mathbb{P}^2,$ then the
scheme $W (\subset \mathbb{P}(T_{\mathbb{P}^2}))$  of zeros of $\phi(s)$ is a rational 
normal surface with exactly one singular point which is of the type $A_1$ (see, \cite{R} \S 4.2)

ii)  If $V$ is a non reduced scheme of length three supported at a point $P$ then 
the scheme $W (\subset \mathbb{P}(T_{\mathbb{P}^2}))$  of zeros of $\phi(s)$ is a rational 
normal surface with exactly one singular point which is of the type $A_2$ (see, \cite{R} \S 4.2)
 
 \end{lem}
 
{\bf Proof}:  From the Lemma \ref{lem1} and its proof we deduce the following:

 If  $V$ is as in i) 
then the scheme $W$ is isomorphic to blow-up of $\PP^2$ at the maximal ideal at
$Q$ and the ideal  $(x,y^2)$ at $P,$ where $x,y$ are  local coordinates at $P.$
Now it can be seen that $W$ is non-singular except one point $P_0$ over $P$
local ring at $P_0$ is isomorphic to localisation of $\CC[x, y, t]/{(xt+y^2)}$ at the maximal
ideal $(x,y,t).$ Thus  $W$ has $A_1$ type singularity at $P_0.$ 

If  $V$ is as in ii) 
then the scheme $W$ is isomorphic to blow-up of $\PP^2$ 
along  the ideal  $(x,y^3)$ at $P,$ where $x,y$ are the local coordinates at $P.$
Hence $W$ is non-singular except one point $P_0$ over $P$
local ring at $P_0$ is isomorphic to localisation of $\CC[x, y, t]/{(xt+y^3)}$ at the maximal
ideal $(x,y,t).$ Thus  $W$ has $A_2$ type singularity at $P_0.$  $\hfill{\Box}$

\vskip3mm

  {\bf Proof of Theorem \ref{thm1}}: Note that  a  hyperplane section of  
  $\mathbb{P}(T_{\mathbb{P}^2})$ in its natural embedding in $\PP^7$ 
  is a two dimensional scheme of degree six and  is the image of the 
  zero scheme of a non zero section 
  of $ {\mathcal O}_{\mathbb{P}(T_{\mathbb{P}^2})}(1).$
   Let $s$ be  a non zero section of 
 $T_{\mathbb{P}^2}$ and $\phi(s)$ be the corresponding section  
 $ {\mathcal O}_{\mathbb{P}(T_{\mathbb{P}^2})}(1).$  
 
 First assume that the zeros of the section 
 $s$  are of the form given by a), b) or c). Then the hypothesis of Lemma \ref{lem1} holds with
 $V_1(s) = \emptyset .$
 Therefore the scheme $W$ of zeros of the section $\phi(s)$ is equal to the blow-up $\mathbb{P}^2$ along
 $V= V_0(s)$ as in Lemma \ref{lem1} and hence irreducible. If $V_0(s)$ as in a) then the $W$ is non-singular
 and is isomorphic to blow-up of  $\PP^2$ along three distinct points. By Lemma \ref{lem3}, if $V$ as in b) then 
 $W$ is rational normal surface with exactly one $A_1$ type singularity  and if $V$ as in c) then 
 $W$ is rational normal surface with exactly one $A_2$ type singularity.
 
 From this we conclude that, if the zeros of the section
 $s$ are of the form given by $a), b)$  or $c) $ in Theorem \ref{thm2}, then 
 the image of the zero scheme of $\phi(s)$
 is an irreducible surface of degree $6$ in $\mathbb{P}^6$ described in the Theorem. 
 
 Next assume that the zeros of the section 
 $s$  are of the form given by $d)$ or $e)$  in Theorem \ref{thm2}. Then the hypothesis of Lemma \ref{lem1} holds with
 $V_1(s) = \ell$ a line in $\mathbb{P}^2$ and $V_0(s) $ is a single point. Since 
 $$T_{\mathbb{P}^2}|_{\ell} = \OOO_{\ell}(1)\oplus  \OOO_{\ell}(2),$$
 from  Lemma \ref{lem1}  we see that the scheme of zeros of the section $\phi(s)$ 
 is equal to $W_0 \cup W_1,$ where $W_0$ is 
 the blow-up $\mathbb{P}^2$ along $V_0(s)$  and $W_1 = \mathbb{P}(\OOO_{\ell}(1)\oplus  \OOO_{\ell}(2)).$ 
 From this we conclude that, if the zeros of the section
 $s$ are of the form given by $d)$ or $e)$ in Theorem \ref{thm2} then  
  the image of the zero scheme of $\phi(s)$
 is a union of two surfaces, namely the images of $W_0$ and $W_1.$
 Since the image of $W_1 = \mathbb{P}(\OOO_{\ell}(1)\oplus  \OOO_{\ell}(2))$
 is a surface of degree 3, we conclude that the image of $W_0$ is also a surface of degree 3, as
 claimed in ii) of the Theorem.  From the Lemma \ref{lem1} it follows that $W_0$ is isomorphic 
 to blow-up of $\PP^2$ at a point. Hence both $W_0$ and $W_1$ are non-singular surfaces.
 It is now easy to see that if the section $s$ is of the form given by $d)$  (respectively, by $e)$) 
 in Theorem \ref{thm2}, then
 the zeros of $\phi(s)$ consists of union of two non-singular irreducible surfaces, intersecting along a
 non-singular rational curve (respectively, intersecting along a union of two non-singular rational curves meeting at a point).
 This completes the proof of the Theorem.  $\hfill{\Box}$
 
\begin{rem} \label{deform}

It is well known that the variety $\mathbb{P}(T_{\mathbb{P}^2})$ can be identified with the homogeneous space ${\text{SL}_3(\mathbb{C})}/B,$
where ${\text{SL}_3}(\mathbb{C})$ is the special linear group of $3\times3$ matrices  with determinant one
over the field of complex numbers and $B$ is the Borel subgroup of upper triangular  matrices. 
Thus $\mathbb{P}(T_{\mathbb{P}^2})$ is 
a homogenous variety. 

1) The nonsingular hyperplane sections  of $\text{SL}_3(\mathbb{C})/B $ described in  the Theorem \ref{thm1} are the well known 
    Del Pezzo Surfaces of degree 6 in $\mathbb{P}^6.$ The Theorem \ref{thm1} describes all the possible 
    degenerations of the Del Pezzo surface of degree 6 inside the homogeneous space $\text{SL}_3(\mathbb{C})/B.$

2)   The very ample linear system that we considered here is given by the sum of two 
homologically inequivalent Schubert divisors in ${\text{SL}_3}(\mathbb{C}).$ These inequivalent Schubert are known as fundamental Schubert divisors in 
${\text{SL}_3(\mathbb{C})}/B.$ In the natural embedding these fundamental Schubert divisors are the degree 3
components of a reducible hyperplane section whose components intersect along the union of two rational curves meeting 
at a single point in the Theorem \ref{thm1}. This can be deduced from  the description of Bruhat order on the symmetric group $S_3$ 
(see. \cite{BM}
Example 1.2.3) and  the description of the  Schubert divisors (see. \cite{BM} Proposition 1.2.1).  
By using the Theorem \ref{thm1} we see that the possible deformations of the 
union of the two fundamental Schubert divisors in  ${\text{SL}_3(\mathbb{C})}/B$ are

 a) union of two irreducible nonsingular
surfaces  $W_1$ and $W_2$ intersecting along a irreducible non-singular rational curve and 
$$W_1 \simeq W_2
\simeq \mathbb{P}(\mathcal{O}_{\mathbb{P}^1} \oplus \mathcal{O}_{\mathbb{P}^1}(1))$$

b) irreducible rational normal surface with exactly one singularity which is either of the type $A_1$ or $A_3.$

c) non-singular Del Pezzo Surfaces of degree 6.

\end{rem}
 
{\bf Acknowledgement: } We would  like to thank University of Lille-1 at Lille, France for its hospitality.
The second named author would  like to thank the University of Artois at Lens, France,
for inviting him  for an academic visit.

 %%%%%%%%%%%%%%%%%%%%%%%%%%%%%%%%%%%%%%%%%%%%%%%%%%%%%%%%%%%%%%%%%%%%%%%%%%

\end{document}